\documentclass[12pt,psamsfonts,leqno,oneside,letterpaper]{amsart}
\usepackage[dvips,text={6.5truein,9truein},left=1truein,top=1truein]{geometry}
\usepackage{amssymb,amsmath,amscd,enumerate}
\usepackage[pdftex]{graphicx}
\usepackage{url}

\usepackage[colorlinks,linkcolor=blue,citecolor=blue,pdfstartview=FitH,bookmarks=false]{hyperref}
\usepackage{amsrefs}
\input xy
\xyoption{all}
\SelectTips{cm}{12}

\usepackage{color}

\theoremstyle{definition}
\newtheorem{para}{}[section]

\newtheorem{remark}[section]{Remark}
\newtheorem{remarks}[section]{Remarks}
\newtheorem{notation}[section]{Notation}
\newtheorem{convention}[section]{Convention}
\newtheorem{definition}[section]{Definition}
\newtheorem{definitions}[section]{Definitions}

\newcommand\Alternatives{\begin{enumerate}[(i)]}
\newcommand\EndAlternatives{\end{enumerate}}
\newcommand\Conditions{\begin{enumerate}[(1)]}
\newcommand\EndConditions{\end{enumerate}}

\theoremstyle{plain}
\newtheorem{theorem}[section]{Theorem}
\newtheorem{lemma}[section]{Lemma}
\newtheorem{proposition}[section]{Proposition}
\newtheorem{corollary}[section]{Corollary}
\newtheorem{conjecture}[section]{Conjecture}
\newtheorem{claim}[equation]{}

\newcommand\Number{\begin{para}}
\newcommand\EndNumber{\end{para}}
\newcommand\Definition{\begin{definition}}
\newcommand\EndDefinition{\end{definition}}
\newcommand\Definitions{\begin{definitions}}
\newcommand\EndDefinitions{\end{definitions}}
\newcommand\Theorem{\begin{theorem}}
\newcommand\EndTheorem{\end{theorem}}
\newcommand\Conjecture{\begin{conjecture}}
\newcommand\EndConjecture{\end{conjecture}}
\newcommand\Remark{\begin{remark}}
\newcommand\EndRemark{\end{remark}}
\newcommand\Remarks{\begin{remarks}}
\newcommand\EndRemarks{\end{remarks}}
\newcommand\Convention{\begin{convention}}
\newcommand\EndConvention{\end{convention}}
\newcommand\Notation{\begin{notation}}
\newcommand\EndNotation{\end{notation}}
\newcommand\Lemma{\begin{lemma}}
\newcommand\EndLemma{\end{lemma}}
\newcommand\Proposition{\begin{proposition}}
\newcommand\EndProposition{\end{proposition}}
\newcommand\Corollary{\begin{corollary}}
\newcommand\EndCorollary{\end{corollary}}
\newcommand\Claim{\begin{claim}}
\newcommand\EndClaim{\end{claim}}
\newcommand\Proof{\begin{proof}}
\newcommand\EndProof{\end{proof}}
\newcommand\Equation{\begin{equation}}
\newcommand\EndEquation{\end{equation}}
\newcommand\NoProof{{\hfill$\square$}}
\newcommand\Bullets{\begin{itemize}}
\newcommand\EndBullets{\end{itemize}}

\newcommand\cala{{\mathcal A}}
\newcommand\calb{{\mathcal B}}
\newcommand\calp{{\mathcal P}}
\newcommand\caln{{\mathcal N}}
\newcommand\arccosh{\mathop{\rm arccosh}}

\newcommand\arcsinh{\mathop{\rm arcsinh}}

\newcommand\ZZ{{\mathbb Z}}
\newcommand\vol{\mathop{\rm Vol}}
\begin{document}

\title{Betti numbers and injectivity radii} 
\author{Marc Culler}
\address{Department of Mathematics, Statistics, and Computer Science
(M/C 249)\\
University of Illinois at Chicago\\
851 S. Morgan St.\\
Chicago, IL 60607-7045}
\email{culler@math.uic.edu}
\thanks{Partially supported by NSF {grants  DMS-0504975 and DMS-0608567}}

\author{Peter B.~Shalen}
\address{Department of Mathematics, Statistics, and Computer Science
(M/C 249)\\
University of Illinois at Chicago\\
851 S. Morgan St.\\
Chicago, IL 60607-7045}
\email{shalen@math.uic.edu}


\maketitle

The theme of this paper is the connection between topological
properties of a closed orientable hyperbolic $3$-manifold $M$ and the
maximal injectivity radius of $M$. In \cite{paradoxical} we showed
that if the first Betti number of $M$ is at least $3$ then the maximal
injectivity radius of $M$ is at least $\log3$. By contrast, the best
known lower bound for the maximal injectivity radius of $M$ with no
topological restriction on $M$ is the lower bound of
$\arcsinh(\frac14)=0.24746\ldots$ due to Przeworski \cite{prez}. One of
the results of this paper, Corollary \ref{STILL no fibroids!}, gives a
lower bound of $0.32798$ for the case where the first Betti number of
$M$ is $2$ and $M$ does not contain a ``fibroid'' (see below).  Our
main result, Theorem \ref{.32whatever}, is somewhat stronger than this

The proofs of our results combine a result due to Andrew Przeworski
\cite{prez} with results from \cite{hyperhakI} and
\cite{hyperhakenII}.

The results of \cite{hyperhakI} and \cite{hyperhakenII} were motivated
by applications to the study of hyperbolic volume, and these
applications were superseded by the results of \cite{ast}. The
applications presented in the present paper do not seem to be
accessible by other methods.

\bigskip
\centerline{$\bf{*}$\qquad$\bf{*}$\qquad$\bf{*}$}
\bigskip
   
As in \cite{hyperhakI}, we define a {\it book of $I$-bundles}
to be a compact, connected, orientable topological 3-manifold (with
boundary) $W$ which has the form $W=\calp\cup \calb$, where
\begin{itemize}
\item $\calp$ is an $I$-bundle over a non-empty compact
2-manifold-with-boundary;
\item each component of $\calb$ is homeomorphic to 
$D^2\times S^1$;
\item the set $\cala =\calp\cap \calb$ is the vertical boundary of the
$I$-bundle $\calp$; and
\item each component of $\cala $ is an  annulus in $\partial \calb$
which is homotopically non-trivial in $\calb$.
\end{itemize}

(Note that this terminology differs slightly from that of \cite{last},
where it is the triple $(W,\calb,\calp)$ that is called a book of $I$-bundles.)

As in \cite{hyperhakI}, we define a {\it fibroid} in a closed,
connected, orientable $3$-manifold $M$ to be a connected
incompressible surface $F$ with the property that each component of
the compact manifold obtained by cutting $M$ along $F$ is a book of
$I$-bundles.  Note that in defining a fibroid to be connected, we are
following the convention of \cite{hyperhakI} rather than that of
\cite{hyperhakenII}.

We define a function $R(x)$ for $x>0$ by
\Equation\label{or whatever your name was}
R(x) = \frac12\arccosh\biggl(
\frac{e^{2x} + 2e^x + 5}{(\cosh \frac{x}{2})(e^x - 1)(e^x +
  3)}\biggr) .
\EndEquation

As in \cite{accs}*{Section 10}, we define a function $V(x)$ for $x>0$
by \Equation\label{this is for you violet} V(x)=\pi x\sinh^2R(x) =
\frac{\pi x}{e^x- 1}\biggl(\frac{e^{2x}+2e^x+5}{2 (e^x+3)\cosh (
  x/2)}\biggr)-\frac{\pi x}2.
\EndEquation

Thus, in a closed hyperbolic $3$-manifold, if a geodesic of length
$\ell$ is the core of an embedded tube of radius $R(\ell)$, then this
tube has volume $V(\ell)$.
 
The following result is implicit in \cite{hyperhakenII}, but we will
supply a proof for the sake of comprehensibility.

\Proposition \label{driven snow} Let $M$ be a closed hyperbolic
$3$-manifold.  Suppose that there is an infinite subset $\caln$ of
$H_2(M;\ZZ)$ such that every element of $\caln$ is represented by some
connected, incompressible surface which is not a fibroid. Let
$\lambda$ be a positive number less than $\log3$. Then either the
maximal injectivity radius of $M$ is at least $\lambda/2$, or $M$
contains a closed geodesic $C$ of length at most $\lambda$ such that
the maximal tube about $C$ has radius at least $R(\lambda)$
and volume at least $V(\lambda)$, where $R(\lambda)$ and
$V(\lambda)$ are defined by (\ref{or whatever your name was}) and
(\ref{this is for you violet}).
\EndProposition

\Proof The hypothesis implies in particular that $H_2(M;\ZZ)$ has
infinitely many primitive elements, and so the first Betti number
$\beta_1(M)$ is at least $2$.  If $\beta_1(M)\ge3$, then according to
\cite{paradoxical}*{Corollary 10.4}, the maximal injectivity radius is
at least $\frac12\log3 > \lambda/2$. We may therefore assume that
$\beta_1(M)=2$.  Hence the quotient of $H_1(M,\ZZ)$ by its torsion
subgroup is a free abelian group $L$ of rank $2$. We let
$h\colon\pi_1(M)\to L$ denote the natural homomorphism.

We distinguish two cases. First consider the case in which $M$
contains a non-trivial closed geodesic $C$ of some length $\ell<\lambda$
such that the conjugacy class represented by $C$ is contained in the
kernel of $h$. Let $T$ denote the maximal embedded tube about
$C$. According to \cite{accs}*{Corollary 10.5} we have $\vol T\ge
V(\lambda)$. If $\rho$ denotes the radius of $T$, this gives
$$ \pi\ell\sinh^2\rho=\vol T\ge
V(\lambda)=\pi\lambda\sinh^2R(\lambda)>\pi\ell\sinh^2R(\lambda),$$
and hence $\rho>R(\lambda)$.

Thus the second alternative of the proposition holds in this case.

We now turn to the case in which no non-trivial closed geodesic of
length $<\lambda$ represents a conjugacy class contained in the
kernel of $h$. Since $M$ is closed, there are only a finite number
$n\ge0$ of conjugacy classes in $\pi_1(M)$ that are represented by
closed geodesics of length $<l$. Let $\gamma_1,\dots,\gamma_n$ be
elements belonging to these $n$ conjugacy classes. Then
$\bar\gamma_i=h(\gamma_i)$ is a non-trivial element of $L$ for
$i=1,\dots,n$. Since $L$ is a free abelian group of rank $2$, there
exists, for each $i\in\{1,\dots,n\}$, a homomorphism $\phi_i$ of 
$L$ onto $\ZZ$ such that $\phi_i(\bar\gamma_i)=0$. Because
$\bar\gamma_i\ne0$, the epimorphism $\phi_i$ is unique up to
sign.

The epimorphism $\phi_i\circ h\colon\pi_1(M)\to\ZZ$ corresponds to a
primitive element of $H^1(M;\ZZ)$, whose Poincar\'e dual in $H_2(M;\ZZ)$
we shall denote by $c_i$. Since the set $\caln\subset H_2(M;\ZZ)$ given
by the hypothesis of the theorem is infinite, there is an element $c$
of $\caln$ which is distinct from $\pm c_i$ for $i=1,\dots,n$. Since
$c\in\caln$ it follows from the hypothesis that there is a connected
incompressible surface $S\subset M$ which represents the homology
class $c$ and is not a fibroid.

We now apply Theorem A of \cite{hyperhakI}, which asserts that if $S$
is a connected non-fibroid incompressible surface in a closed,
orientable hyperbolic $3$-manifold $M$, and if $\lambda$ is any
positive number, then either (i) $M$ contains a non-trivial closed
geodesic of length $<\lambda$ which is homotopic in $M$ to a closed
curve in $M-S$, or (ii) $M$ contains a hyperbolic ball of radius
$\lambda/2$.  In the present situation, with $\lambda$ chosen as
above, we claim that alternative (i) of the conclusion of Theorem A of
\cite{hyperhakI} cannot hold.

Indeed, suppose that $C$ is a non-trivial closed geodesic of
length $<\lambda$ with the properties stated in (i). Since $C$ has
length $<\lambda$, the conjugacy class represented by $C$ contains
$\gamma_i$ for some $i\in\{1,\dots,n\}$. Since $C$ is homotopic to
a closed curve in $M-S$ it follows that the image of $\gamma_i$ in
$H_1(M;\ZZ)$ has homological intersection number $0$ with $c$.
Thus if $\psi\colon\pi_1(M)\to\ZZ$ is the homomorphism
corresponding to the Poincar\'e dual of $c$, we have
$\psi(\gamma_i)=0$. Now since $L$ is the quotient of $H_1(M)$ by
its torsion subgroup, $\psi$ factors as $\phi\circ h$, where $\phi$
is some homomorphism from $L$ to $\ZZ$. Since $c$ is primitive,
$\psi$ is surjective, and hence so is $\phi$. But we have
$\phi(\bar\gamma_i)=\psi(\gamma_i)=0$. In view of the uniqueness
that we observed above for $\phi_i$, it follows that
$\phi=\pm\phi_i$, so that $\psi=\pm\phi_i\circ h$ and hence $c=\pm
c_i$. This contradicts our choice of $c$.

Hence (ii) must hold. This means that the maximal injectivity radius
of $M$ is at least $\lambda/2$. Thus the first alternative of the
proposition holds in this case.
\EndProof

\Proposition \label{pure .3279}
Let $M$ be a closed hyperbolic $3$-manifold.  Suppose that there is an
infinite subset $\caln$ of $H_2(M;\ZZ)$ such that every element of
$\caln$ is represented by some connected, incompressible surface which
is not a fibroid. Then the maximal injectivity radius of $M$ exceeds
$0.32798$.
\EndProposition

\Proof
We set $\lambda=2\times0.32798=0.65596$. 

According to Proposition \ref{driven snow}, either the maximal
injectivity radius of $M$ is at least $\lambda/2$ --- so that the
conclusion of the theorem holds --- or $M$ contains a closed geodesic
$C$ of length at most $\lambda$ such that the maximal tube about $C$
has volume at least $V(\lambda)$ , where $V(\lambda)$ is defined by
(\ref{this is for you violet}).
In the latter case, if  $R$ denotes the radius of $T$, we have
$$R \ge R(\lambda)= 0.806787\ldots .$$

Now according to \cite{prez}*{Proposition 4.1}, the maximal injectivity
radius of $M$ is bounded below by
$$
\arcsinh\bigg(\frac{\tanh R}2\bigg)>
\arcsinh\bigg(\frac{\tanh0.806787 }2\bigg)
>0.32799.$$ 
This gives the conclusion of the theorem in this case.
\EndProof

\Theorem \label{.32whatever}
Let $M$ be a closed hyperbolic $3$-manifold.  Suppose that there is an
infinite set $\caln$ of primitive elements of $H_2(M;\ZZ)$ such that
no element of $\caln$ is represented by a (connected) fibroid.  Then
the maximal injectivity radius of $M$ exceeds $0.32798$.
\EndTheorem

\Proof If $\pi_1(M)$ has a non-abelian free quotient, then by
\cite{hyperhakenII}*{Theorem 1.3}, the maximal injectivity radius of
$M$ is at least  $\frac12\log3 = 0.549\ldots$. Now suppose
that $\pi_1(M)$ has no non-abelian free quotient.  If $\caln$ is the
set given by the hypothesis of Theorem \ref{.32whatever}, it now
follows from \cite{hyperhakenII}*{Proposition 2.1} that every element
of $\caln $ is represented by a connected incompressible surface,
which by hypothesis cannot be a fibroid. Thus $\caln$ has the
properties stated in the hypothesis of Proposition \ref{pure
  .3279}. The latter result therefore implies that the maximal
injectivity radius of $M$ exceeds $0.32798$.
\EndProof

If $M$ is a $3$-manifold whose first Betti number is at least $2$,
then $H_2(M;\ZZ)$ has infinitely many primitive elements. If a
non-trivial element of $H_2(M;\ZZ)$ is represented by a connected
surface it must be primitive, since it has intersection number $1$
with a class in $H_1(M;\ZZ)$.  Hence Theorem \ref{.32whatever} implies:

\Corollary \label{STILL no fibroids!}  
Let $M$ be a closed hyperbolic $3$-manifold.  Suppose that the first
 Betti  number of $M$ is at least $2$, and that $M$ contains
no non-separating fibroid.  Then the maximal injectivity radius of $M$
exceeds $0.32798$.
\EndCorollary
\NoProof

\end{document}